\documentclass{article}

\usepackage{amsthm,amssymb}
\catcode`\@=11
\@addtoreset{equation}{section}

\newtheorem{thm}{Theorem}[section]
\newtheorem{prop}{Proposition}[section]

\newtheorem{rem}{Remark}[section]
\newtheorem{lem}{Lemma}[section]

\newcommand{\eps}{\varepsilon}
\newcommand{\loc}{{\mathrm{loc}}}
\newcommand{\sgn}{{\mathrm{sgn}\,}}

\newcommand{\R}{{\mathbb{R}}}
\newcommand{\kbf}{{\mathbf{k}}}

\newcommand{\xbf}{{\mathbf{x}}}

\newcommand{\tbf}{{\mathbf{t}}}
\newcommand{\sbf}{{\mathbf{s}}}
\newcommand{\qbf}{{\mathbf{q}}}
\newcommand{\C}{{\mathbb{C}}}
\newcommand{\Z}{{\mathbb{Z}}}
\newcommand{\T}{{\mathbb{T}}}
\newcommand{\N}{{\mathbb{N}}}
\newcommand{\Q}{{\mathbb{Q}}}
\newcommand{\calA}{\mathcal{A}}

\newcommand{\calB}{\mathcal{B}}
\newcommand{\calJ}{\mathcal{J}}
\newcommand{\calH}{\mathcal{H}}
\newcommand{\calL}{\mathcal{L}}
\newcommand{\calM}{\mathcal{M}}
\newcommand{\calN}{\mathcal{N}}

\newcommand{\calP}{\mathcal{P}}
\newcommand{\calQ}{\mathcal{Q}}

\newcommand{\Real}{{\mathrm{Re}\,}}

\title{Variational approach to second species periodic solutions of Poincar\'e of the 3 body problem
\author{Sergey Bolotin\thanks{Supported by the Programme
``Mathematical Control Theory'' of RAS and  RFBR
under grant \# 08-01-00681-a.
The first draft of this paper was written while the first author was visiting CIRM in November 2008.}\\
Department of Mathematics\\
University of Wisconsin--Madison\\
and\\
Moscow Steklov Mathematical Institute
\and Piero Negrini\\
Department of Matematics \\
Sapienza, University of Rome} }

\begin{document}

\maketitle

\centerline{Dedicated to Ernesto Lacomba on the occasion of his 65th birthday}

\begin{abstract}
We consider the plane 3 body problem with 2 of the masses small.
Periodic solutions with near collisions of small bodies were named
by Poincar\'e second species periodic solutions.  Such solutions
shadow chains of collision orbits of 2 uncoupled Kepler problems.
Poincar\'e only sketched the proof of the existence of second
species solutions. Rigorous proofs appeared much later and only for
the restricted 3 body problem. We develop a variational approach to
the existence of second
species periodic solutions for the nonrestricted 3 body problem.
As an application, we give a rigorous proof of the existence
of a class of second species solutions.
\end{abstract}

\section{Introduction}

Consider the plane 3-body problem   with masses $m_1,m_2,m_3$.
Suppose that $m_3$ is much larger than $m_1,m_2$, i.e.\
$\mu=(m_1+m_2)/m_3$ is a small parameter. Then
$$
m_1/m_3=\mu\alpha_1,\quad m_2/m_3=\mu\alpha_2,\quad \alpha_1+\alpha_2=1.
$$
Let $q_1,q_2\in\R^2$ be positions of   $m_1,m_2$ with respect to
$m_3$ (Poincar\'e's heliocentric coordinates) and $p_1,p_2\in\R^2$
their scaled momenta.  The motion of  $m_1,m_2$ with respect to
$m_3$ is described by a Hamiltonian system ($H_\mu$) with
Hamiltonian
\begin{equation}
\label{eq:Hmu}
H_\mu(q,p)=H_0(q,p)+\mu \left(
\frac{|p_1+p_2|^2}2 - \frac{\alpha_1\alpha_2}{|q_1-q_2|}\right),
\end{equation}
where
$$
H_0(q,p)=   \frac{|p_1|^2}{2 \alpha_1}+\frac{|p_2|^2}{2 \alpha_2} -\frac{
\alpha_1}{|q_1|}  -\frac{ \alpha_2}{|q_2|}.
$$

The Hamiltonian $H_\mu$  is a small perturbation of the Hamiltonian
$H_0=H_1+H_2$   describing 2 uncoupled Kepler problems. The
configuration space of system $(H_0)$ is
$U^2=U\times U$, where $U=\R^2\setminus\{0\}$.
The configuration space of the perturbed system $(H_\mu)$ is
$U^2\setminus \Delta$, where
$$
\Delta=\{q=(q_1,q_2)\in U^2:q_1=q_2\}
$$
represents collisions of $m_1,m_2$.

System $(H_\mu)$ has energy integral and the integral of angular
momentum\footnote{We identify $\R^2$ with $\C$, so multiplication by $i$ is rotation by $\pi/2$.}
$$
G(q,p)=G_1+G_2=iq_1\cdot p_1+iq_2\cdot p_2
$$
corresponding to the rotational symmetry $e^{i\theta}:\R^2\to \R^2$.

System $(H_0)$ has additional first integrals $H_1, H_2$ and
$G_1,G_2$
 -- energies and angular momenta of  $m_1,m_2$. In the domain
$$
P=\{(q,p)\in U^2\times \R^4: H_1,H_2<0,\; G_1,G_2\ne 0\},
$$
orbits of  $m_1,m_2$ are Kepler ellipses, and
solutions of system ($H_0$) are quasiperiodic with 2 frequencies.

Let  $R\subset P$  be the regular domain -- the set of points
in $P$ such that the corresponding Kepler ellipses do not cross.
For small $\mu>0$ solutions of system ($H_\mu$) in $R$ are
$O(\mu)$-approximated by solutions of system ($H_0$)  on finite time
intervals independent of $\mu$. By the classical perturbation
theory, away from resonances the same is true on longer time
intervals. Moreover, as proved by Arnold \cite{Arnold}, for small
$\mu>0$,
 system ($H_\mu$) has  a large measure of invariant 2-dimensional KAM tori
 on which solutions are quasiperiodic, and thus well approximated
 (modulo rotation) by solutions of system $(H_0)$ on infinite time intervals.

In the singular domain $S\subset P$, where the corresponding
Kepler ellipses cross, the classical perturbation theory does not
work. Indeed, for almost any initial condition in $S$, solution
of system $(H_0)$ is quasiperiodic with incommensurable frequencies,
and so eventually $m_1,m_2$ simultaneously approach an intersection
point of Kepler ellipses. Then the perturbation in
(\ref{eq:Hmu}) becomes large, and so it can not be ignored even in
the first approximation in $\mu$.

For small $\mu>0$  solutions of the 3 body problem $(H_\mu)$ in $S$
can be described as follows. The  bodies $m_1,m_2$ move along nearly
Kepler ellipses and after many revolutions they almost collide. Then
they start moving  along a new pair  of nearly Kepler orbits. If the
new energies $H_1,H_2$ are both negative, so the new Kepler orbits
are ellipses, then $m_1,m_2$ will again nearly collide after many
revolutions, and the process repeats itself. Thus  almost collision
solutions of the 3-body problem $(H_\mu)$ shadow chains of collision
orbits of system $(H_0)$.

  Almost collision periodic solutions of system
$(H_\mu)$ were first studied by Poincar\'e  in New Methods of Celestial Mechanics.
Poincar\'e named them second species periodic solutions. However, he did not
provide a rigorous existence proof. Rigorous proofs appeared much
later (see e.g.\ \cite{Mar-Nid,Bol-Mac,Bol:ELL}) and only for the
restricted 3 body problem, circular  and elliptic. In
\cite{Bol-Mac,Bol:Nonlin} also chaotic second species solutions of
the circular and elliptic restricted problem were studied.

The goal of this paper is to develop a variational approach to
almost collision  periodic orbits of the {\it nonrestricted} 3 body
problem. As an application, we will give a rigorous proof of the
existence of a class of almost collision periodic orbits. Chaotic
almost collision orbits will be studied in another paper.

\begin{rem}
It is possible to fix the value of angular momentum $G$ and reduce
rotational symmetry. Then we obtain a Hamiltonian system with 3
degrees of freedom.  However,
since reduction of the rotational symmetry considerably complicates
the Hamiltonian, it is simpler to work with the original Hamiltonian
system $(H_\mu)$ with 4 degrees of freedom.
\end{rem}

\begin{rem}
We consider only near collisions of small masses $m_1$ and $m_2$ and
exclude near collisions of $m_1,m_2$ with $m_3$. In particular,
triple collisions are excluded. It is well known that double
collisions can be regularized, but the Levi-Civita regularization
becomes singular as $\mu\to 0$. Understanding this singularity is
the base for our methods. Levi-Civita regularization was previously
used to study second species solutions for the restricted 3 body
problem, see e.g.\ \cite{Mar-Nid,Bol-Mac,FNS}.
\end{rem}

\begin{rem}
Main results of this paper
hold for more general Hamiltonians with singularity, for example
$$
H_\mu(q,p)=\frac{|p_1|^2}{2 a_1(q)}+\frac{|p_2|^2}{2 a_2(q)}
+g(q,\mu)   - \frac{\mu f(q,\mu)}{|q_1-q_2|},\qquad a_1,a_2,f>0,
$$
where all functions are smooth without singularity at $q_1=q_2$.
\end{rem}

\section{Main results}

A solution   of the Hamiltonian system $(H_\mu)$ is determined by
its projection  to the configuration space $U^2\setminus\Delta$
which will be called a trajectory.  Let $L_\mu$ be the Lagrangian
corresponding to $H_\mu$. A $T$-periodic trajectory $\gamma:\R\to
U^2\setminus\Delta$ is a critical point of the Hamilton action
functional
\begin{equation}
\label{eq:Amu}
A_\mu(T,\gamma)=\int_0^TL_\mu(\gamma(t),\dot\gamma(t))\,dt
\end{equation}
on the space of $T$-periodic $W^{1,2}_\loc$ curves in
$U^2\setminus\Delta$. We write $T$ explicitly in $A_\mu(T,\gamma)$
because later it will become a variable.

Any trajectory of system $(H_0)$ has the form
$\gamma=(\gamma_1,\gamma_2)$, where $\gamma_j$ is a trajectory of
the Kepler problem. For $\mu=0$ the action functional $A=A_0$
splits:
\begin{equation}
\label{eq:action}
A(T,\gamma)= \alpha_1B(T,\gamma_1)+\alpha_2 B(T,\gamma_2),\qquad \gamma=(\gamma_1,\gamma_2),
\end{equation}
where
\begin{equation}\label{eq:B}
B(T,\sigma)=\int_0^T \left(\frac{|\dot
\sigma(t)|^2}{2}+\frac{1}{|\sigma(t)|}\right)dt
\end{equation}
is the action functional of the Kepler problem on the space
of $T$-periodic $W^{1,2}_\loc$ curves $\sigma:\R\to U$.

We have to consider also trajectories of system $(H_0)$  with collisions. A $T$-periodic curve
$\gamma=(\gamma_1,\gamma_2):\R\to  U^2$ is called a periodic
$n$-{\it collision chain} if there exist time moments
\begin{equation}
\label{eq:tbf}
\tbf=(t_1,\dots,t_n),  \qquad t_1<\dots<t_n<t_{n+1}=t_1+T,
\end{equation}
such that:
\begin{itemize}
\item
$\gamma$ has collisions at $t=t_j$, so that
$$
\gamma(t_j)=(x_j,x_j)\in \Delta,\qquad \gamma_1(t_j)=\gamma_2(t_j)=x_j.
$$
\item
$\gamma|_{[t_j,t_{j+1}]}$   is a trajectory of system $(H_0)$ which
will be called a {\it collision orbit.}
\item
Momentum $p(t)=(p_1(t),p_2(t))=(\alpha_1\dot\gamma_1(t),\alpha_2\dot\gamma_2(t))$
changes    at collisions so that the total momentum $y=p_1+p_2$ is continuous:
\begin{equation}
\label{eq:mom}
y(t_j+ 0)=y(t_j-0)=y_j.
\end{equation}
Equivalently, the jump of momentum $p(t_j+0)-p(t_j-0)$ is orthogonal to $\Delta$ at
$\gamma(t_j)$.
\item The total energy
$H_0$  does not change at collision:
\begin{equation}
\label{eq:en}
H_0(\gamma(t_j),p(t_j+0))=
H_0(\gamma(t_j),p(t_j-0)).
\end{equation}
\end{itemize}

By (\ref{eq:mom}), the total angular momentum $G=G_1+G_2$ is
preserved at collisions:
$$
G(\gamma(t_j),p(t_j\pm 0))=ix_j\cdot y_j.
$$
Hence  $G$ is constant along $\gamma$.
By (\ref{eq:en}), the total energy $H_0=H_1+H_2=E$ is also constant
along $\gamma$, but not the energies $H_1,H_2$ of $m_1,m_2$, or
their angular momenta $G_1,G_2$.

A  collision chain $\gamma$ is a broken trajectory of system $(H_0)$
-- a concatenation of collision orbits with reflections from
$\Delta$. However, unlike for ordinary billiard systems, $\Delta$
has codimension 2 in the configuration space $U^2$, so the change of
the normal component of the momentum at collision is not uniquely
determined. Thus there is no direct interpretation of collision
chains as trajectories of a dynamical system. Such an interpretation
is given later on.

Collision chains are limits of trajectories of system ($H_\mu$) which approach
collisions as $\mu\to 0$. Indeed, we have:

\begin{prop}
\label{prop:limit} Let $\gamma_\mu$ be a $T_\mu$-periodic trajectory of
system $(H_\mu)$  which uniformly converges, as $\mu\to 0$, to a
$T$-periodic curve $\gamma$. If $\gamma([0,T])\cap\Delta$ is a finite set,
then $\gamma$ is a periodic collision chain.
\end{prop}

We say that $\gamma_\mu$ is an almost collision orbit shadowing the
collision chain $\gamma$.
A similar statement holds for  nonperiodic collision chains.

Intuitively, Proposition \ref{prop:limit}
is almost evident: a near collision of $m_1,m_2$  lasts a short time
during which the influence of the non-colliding body $m_3$ is
negligible. Then   $m_1,m_2$ form a 2 body problem, so their total
momentum $y=p_1+p_2$ and total energy $H_0=H_1+H_2$ are almost
preserved. This yields  (\ref{eq:mom})--(\ref{eq:en}).
 This can be made into a rigorous proof, see e.g.\
\cite{Alexeyev}. A better way to prove Proposition \ref{prop:limit}
is by using the Levi-Civita regularization, see section
\ref{sec:Levi}. \qed

\medskip

Collision chains can be characterized as extremals of Hamilton's
action functional (\ref{eq:action}).  Let $\Omega_n^T$ be the set of
$\omega=(\tbf,T,\gamma)$, where $\tbf=(t_1,\dots,t_n)$ satisfies
(\ref{eq:tbf}) and  $\gamma:\R\to U^2$  is a $T$-periodic
$W^{1,2}_\loc$ curve such that $\gamma(t_j)=(x_j,x_j)\in \Delta$.
The collision times $t_j$ and collision points $x_j$  are not fixed.
Then $\Omega_n^T$ can be identified with an open set in a  Hilbert
space (see (\ref{eq:Omega}) and collision times $t_j$ and collision
points $x_j$ are smooth functions on $\Omega_n^T$.

\begin{rem} If  $\gamma(t)\notin\Delta$ for $t\ne t_j$,
then time moments $t_j$ are determined by the curve $\gamma$, i.e.\
the projection $\Omega_n^T\to W^{1,2}(\R/T\Z,U^2)$,
$(\tbf,T,\gamma)\to\gamma$, is injective. Then
$\omega=(\tbf,T,\gamma)$  is determined by $\gamma$. But $t_j$ are
not continuous functions of $\gamma$, so we have to include the
variables $\tbf$ in the definition of $\Omega_n^T$.
\end{rem}

\begin{rem} In the one-dimensional calculus of variations Hilbert  spaces
are unnecessary: at least locally function spaces can be replaced by
finite dimensional subspaces of broken extremals. We will use this approximation
later on. However, in this section we use conventional $W^{1,2}$ setting.
\end{rem}

The action functional
$\calA(\omega)=A(T,\gamma)$ is a smooth function on $\Omega_n^T$.

\begin{prop}\label{prop:Ham}
$\gamma$ is a $T$-periodic $n$-collision chain iff $\omega=(\tbf,T,\gamma)$
is a critical point of the Hamilton action $\calA$ on $\Omega_n^T$.
\end{prop}

\proof If $\omega=(\tbf,T,\gamma)$ is a critical point of $\calA$ on $\Omega_n^T$,
then each segment $\gamma|_{[t_{i-1},t_j]}$ is a collision orbit of
system ($H_0$). Then by the first variation formula \cite{AKN},
$$
d \calA(\omega)=\sum_{j=1}^n(\Delta h_j\,d t_j-\Delta y_j\cdot d x_j),
$$
where $\Delta h_j$ and $\Delta y_j$ are jumps of energy and  total momentum:
$$
\Delta y_j=y_j^+-y_j^-,\quad y_j^\pm=y(t_j\pm 0), \quad \Delta
h_j=h_j^+-h_j^-,\quad  h_j^\pm=H_0(\gamma(t_j),p(t_j\pm 0)).
$$
Since the differentials $dt_j,dx_j$ are independent, critical points of $\calA$
satisfy $\Delta h_j=0$ and $\Delta
y_j=0$. This implies (\ref{eq:mom})--(\ref{eq:en}).
Converse is also evident.
\qed

\medskip

For collision chains with fixed energy $H_0=E$ we use  Maupertuis's
variational principle \cite{AKN}. Let
$\Omega_n=\cup_{T>0}\Omega_n^T$. The period $T>0$  is now a smooth
function on  $\Omega_n$.  Hamilton's action is replaced by the
Maupertuis action functional
\begin{equation}
\label{eq:AE} \calA^E(\omega)=A^E(T,\gamma)=A(T,\gamma)+ET,\qquad \omega=(\tbf,T,\gamma).
\end{equation}
If $\gamma$ is a collision chain with energy $E$, then
$$
A^E(T,\gamma)=\int_\gamma p\cdot dq
$$
is the classical Maupertuis action.
The  Maupertuis principle  for collision chains is as follows:

\begin{prop}\label{prop:Maup}
$\gamma$ is a $T$-periodic $n$-collision chain of energy $E$ iff
$\omega=(\tbf,T,\gamma)$ is an extremal  of the Maupertuis functional $\calA^E$
on $\Omega_n$.
\end{prop}

Due to time and rotation invariance, critical points of the action
functional are degenerate. To eliminate time degeneracy, we identify
collision chains which differ by time translation
$$
\gamma(t)\to\gamma(t-\tau), \quad t_j\to t_j+\tau,\qquad \tau\in\R.
$$
This defines a group action of $\R$ on $\Omega_n$. The
coordinates on the quotient space $\widetilde\Omega_n=\Omega_n/\R$
can be defined as follows. Let
$$
s_j=t_{j+1}-t_j>0,\quad \sigma_j(\tau)=\gamma_j(t_j+\tau s_j),\qquad
0\le \tau\le 1.
$$
Denote
$$
\sbf=(s_1,\dots,s_n)\in\R_+^n,\qquad
\sigma=(\sigma_1,\dots,\sigma_n)\in W_n,
$$
where
$$
W_n=\{\sigma=(\sigma_1,\dots,\sigma_n):\sigma_j\in
W^{1,2}([0,1],U^2),\quad \sigma_j(1)=\sigma_{j+1}(0)\in\Delta\}.
$$
The map map $\Omega_n\to\R_+^n\times W_n$, $(\tbf,T,\gamma)\to
(\sbf,\sigma)$, makes it possible to identify  $\widetilde\Omega_n$
with $\R_+^n\times W_n$. We can represent $\sigma_j$ by
$$
\hat\sigma_j\in W_0^{1,2}([0,1],U^2),\qquad \hat\sigma_j(\tau)=\sigma_j(\tau)-(1-\tau)x_j-\tau x_{j+1}.
$$
Then $\sigma$ is determined  by $(\xbf,\hat\sigma)$, where
$\xbf=(x_1,\dots,x_n)\in U^n$ and
$\hat\sigma=(\hat\sigma_1,\dots,\hat\sigma_n)$. Hence $W_n\cong
U^n\times W_0^{1,2}([0,1],U^{2n})$. Finally
\begin{equation}
\label{eq:Omega} \widetilde\Omega_n \cong \R_+^n\times W_n\cong
\R_+^n\times U^n\times W_0^{1,2}([0,1],U^{2n}),\quad \Omega_n
\cong\widetilde\Omega_n\times\R.
\end{equation}

The action functional  gives a smooth function $\calA(\sbf,\sigma)$
on $\widetilde\Omega_n$, invariant under rotations:
$\calA(\sbf,e^{i\theta}\sigma)=\calA(\sbf,\sigma)$.
We deal with the rotation degeneracy later on.

Often it is convenient to use parametrization independent Jacobi's form
of the Maupertuis action functional -- the length of $\gamma$ in the Jacobi's metric $ds_E$:
\begin{eqnarray*}
J^E(\gamma)=\int_\gamma ds_E,\qquad
ds_E=\max_p\{p\cdot dq: H_0(q,p)=E\}\\
=\sqrt{2(E+\alpha_1|q_1|^{-1}+\alpha_2|q_2|^{-1})(\alpha_1|dq_1|^2+\alpha_2|dq_2|^2)}.
\end{eqnarray*}
Then $J^E(\gamma)\le A^E(T,\gamma)$ and if $\gamma$ is parametrized
so that $H_0\equiv E$, then $A^E(T,\gamma)=J^E(\gamma)$. Thus up to
parametrization, extremals of $J^E$ and $A^E$ are the same. For a
collision chain $\gamma$ corresponding to
$(\sbf,\sigma)\in\widetilde\Omega_n$,
$$
J^E(\gamma)=\calJ^E(\sigma)=\sum_{j=1}^nJ^E(\sigma)
$$
is a function on $W_n$. We obtain

\begin{prop}\label{prop:Jacobi}
If $\gamma$ is a periodic $n$-collision chain of energy $E$ then the
corresponding $\sigma\in W_n$ is an extremal  of the Jacobi action
$\calJ^E$. If $\sigma$ is an extremal of $\calJ^E$, then if each
$\sigma_j$ is reparametrized so that $H_0\equiv E$, the
corresponding $\gamma$ is  a periodic $n$-collision chain.
\end{prop}

We will not use Jacobi's variational principle since $\calJ^E$ is not a smooth function.
A discrete version of Jacobi's action is smooth and we will use it later on.

Similar variational principles hold for collision chains periodic in
a rotating coordinate frame: $\gamma(t+T)=e^{i\Phi}\gamma(t)$ for some quasiperiod $T$
and phase $\Phi$. We call $\gamma$ {\it periodic modulo rotation}.
If $\Phi\notin 2\pi\Q$, then $\gamma$ is
quasiperiodic in a fixed coordinate frame.

The corresponding function space is defined as follows.
Let
$\hat\Omega_n$
be the set of all $(\tbf,T,\gamma,\Phi)$, where $\tbf,T$ are as before, $\Phi\in\R$ and
the $W_\loc^{1,2}$ curve $\gamma:\R\to U^2$ satisfies
$\gamma(t+T)=e^{i\Phi}\gamma(t)$ and $\gamma(t_j)\in\Delta$.
Then $\hat\Omega_n$ can be identified with an open set in a Hilbert space
and $t_j,x_j,T,\Phi$ are smooth functions on $\hat\Omega_n$. In fact $\hat\Omega_n\cong\Omega_n\times\R$.
Indeed, $\hat\gamma(t)=e^{-it\Phi/T}\gamma(t)$ is a $T$-periodic curve, so $(\tbf,T,\hat\gamma)\in\Omega_n$.

Define the Maupertuis--Routh action functional on $\hat\Omega_n$ by
\begin{equation}
\label{eq:AEG} \calA^{EG}(\tbf,T,\gamma,\Phi)=A(T,\gamma)+ET-G\Phi.
\end{equation}
This is a smooth function on $\hat\Omega_n$ and we have:

\begin{prop} $\gamma$ is a periodic modulo rotation collision chain
with energy $E$ and angular momentum $G$ iff $(\tbf,T,\gamma,\Phi)$ is a
critical point of the functional $\calA^{EG}$ on $\hat\Omega_n$.
\end{prop}

\begin{rem}
It seems natural to take $\Phi\in\T=\R/2\pi\Z$ since $\Phi+2\pi$ gives the same collision chain.
But then  the functional $\calA^{EG}$ will be multivalued: defined modulo $2\pi G$.
\end{rem}

\begin{rem}
The name of the functional $\calA^{EG}$ is motivated as follows.
One can perform Routh's reduction \cite{AKN} of the rotational symmetry for fixed $G$
replacing the configuration space $U^2$ by $\widetilde U^2=U^2/\T\cong \R_+^2\times\T$
and the Lagrangian by the so called Routh function.
Then the functional $\calA^{EG}$ becomes the Maupertuis
functional for the reduced Routh system. Probably this observation is due to Birkhoff \cite{Birkhoff}.
However, Routh's reduction makes the Lagrangian more complicated,
so we do not use it in this paper.
\end{rem}

There are  several other possible variational principles for collision chains (for example,
we may fix the phase $\Phi$), but in the present paper we will use
only the ones given above.

Sufficient condition for the existence of a periodic  orbit of
system $(H_\mu)$ shadowing a given collision chain  $\gamma$
requires that the chain is nontrivial in the following sense.
Let $u(t)= \gamma_2(t)-\gamma_1(t)$ and let
\begin{equation}
\label{eq:v} v(t)=\alpha\dot u(t)=\alpha_1p_2(t)-\alpha_2p_1(t),\qquad  \alpha=\alpha_1\alpha_2,
\end{equation}
be the scaled relative velocity of $m_1,m_2$. Let $v_j^\pm=v(t_j\pm 0)$
be relative collision velocities. Since
\begin{equation}
\label{eq:v_i}
h_j^\pm=|y_j|^2/2+|v_j^\pm|^2/2\alpha
-|x_j|^{-1}=E,
\end{equation}
equations (\ref{eq:mom})--(\ref{eq:en}) imply that
$|v_j^+|=|v_j^-|$: relative speed is preserved at collision.
We impose two essentially equivalent conditions:

\medskip

\noindent {\bf Direction change condition.} Relative collision velocity
changes direction at collision: $\Delta v_j=v_j^+- v_j^-\ne 0$,
$j=1,\dots,n$. In particular, $v_j^\pm\ne 0$.

\medskip

\noindent {\bf No early collisions condition.}  $\gamma(t)\notin\Delta$ for
$t\ne t_j$.

\medskip

If the direction change condition is not satisfied at some $t_j$,
then $\dot\gamma(t_j-0)=\dot\gamma(t_j+0)$, and so
$\gamma|_{[t_{j-1},t_{j+1}]}$ is a smooth trajectory  of system
$(H_0)$. Deleting the collision time  moment $t_j$ we obtain a $(n-1)$-collision
chain violating no early collisions condition.

Conversely, if $\gamma$ is a $n$-collision chain violating no early collisions condition,
then adding an extra collision time moment,
we obtain a $(n+1)$-collision chain violating the changing direction condition.
From now on we add these two equivalent conditions to the definition of a collision chain.

\begin{rem}
The changing direction condition implies that almost collision
orbits $\gamma_\mu$ shadowing the collision chain $\gamma$ come
$O(\mu)$-close to collision. Often almost collision orbits discussed
in Astronomy come close to collision, but not too close, for example
$O(\mu^\nu)$-close with $\nu\in(0,1)$, see e.g.\
\cite{Gomez,Perko,FNS}. Such  orbits change direction at near
collision, but this change  is small as $\mu\to 0$. Then the
corresponding collision chains  do not satisfy the changing
direction condition: $ \Delta v_j=v_j^+-v_j^-=0$. Our methods do not
work for such almost collision orbits.
\end{rem}

The changing direction condition makes it possible to construct a
shadowing orbit $\gamma_\mu$ of system $(H_\mu)$, but it does not
prevent  $\gamma_\mu$ from  having regularizable double collisions
of $m_1,m_2$. To exclude such collisions we need to impose an extra
condition:

\medskip

\noindent{\bf No return condition.} $v_j^++v_j^-\ne 0$, $j=1,\dots,n$.

\medskip

But this condition is not as crucial as the changing direction
condition, so we do not include it in the definition of a collision
chain.

To construct shadowing orbits we also need some nondegeneracy
assumptions. We say that a $T$-periodic $n$-collision chain $\gamma$
with energy $E$ is {\it nondegenerate} if
$\omega=(\tbf,T,\gamma)\in\Omega_n$ is a {\it nondegenerate modulo
symmetry} critical point of the Maupertuis action $\calA^E$ on
$\Omega_n$.

Due to time translation and rotation invariance, critical points of
$\calA^E$ are all degenerate: the group action $\gamma(t)\to e^{i\theta}\gamma(t-\tau)$ of
$\R\times\T$ preserves $\calA^E$. We say that
$\omega=(\tbf,T,\gamma)\in\Omega_n$ is nondegenerate modulo symmetry
if the nullity of the quadric form $d^2\calA(\omega)$  on
$T_{\omega}\Omega_n$ is 2 -- the lowest possible. Equivalently, the
manifold $M\subset\Omega_n$ obtained from $\omega$ by the action of
the group $\R\times\T$ is  a nondegenerate critical manifold.
Nondegeneracy modulo symmetry is equivalent to nondegeneracy of the
corresponding critical point on
$\Omega_n/(\R\times\T)\cong\widetilde\Omega_n/\T$.

As usual in the classical calculus of variations, the Hessian
operator corresponding to $d^2\calA^E(\omega)$ is a sum of
invertible and compact operators on the Hilbert space
$T_{\omega}\Omega_n$, so nondegeneracy modulo symmetry implies that the Hessian has
bounded inverse on $T_{\omega}\Omega_n/T_{\omega}M$. In
fact, at least locally, $\calA^E$ can be reduced to a finite dimensional
discrete action functional (see section  \ref{sec:discr}), so all
Hilbert spaces involved are essentially finite dimensional.

Now two main results will be formulated.

\begin{thm}\label{thm:E1}
Let $\gamma$ be a nondegenerate $T$-periodic collision chain
with energy $E$.  Then for small $\mu>0$ there exists a
$T_\mu$-periodic orbit $\gamma_\mu$ of system $(H_\mu)$ with energy
$E$ which $O(\mu)$ shadows $\gamma$:
$$
T_\mu=T+O(\mu),\qquad \gamma_\mu(t)=\gamma(t)+O(\mu),\qquad t\in[0,T].
$$
If $\gamma$ satisfies the no return condition, then $\gamma_\mu$ has
no collisions and there exist $0<a<b$ independent of $\mu$ such that
\begin{equation}
\label{eq:O(mu)}
\mu a\le d(\gamma_\mu(t_j),\Delta)\le \mu b.
\end{equation}
\end{thm}

Due to time translation and rotation symmetry, 4 multipliers of $\gamma_\mu$ (eigenvalues of the
linear symplectic Poincar\'e map of $\R^8$) are equal to 1.
Nontrivial multipliers are $\lambda_1, \lambda_1^{-1},
\lambda_2,\lambda_2^{-1}$,   where $\lambda_1(\mu)$ is real and large of order
$|\ln\mu|$, and $\lambda_2(\mu)$ has a limit as $\mu\to 0$,
where $\lambda_2(0)\ne 1$ is complex with $|\lambda_2(0)|=1$ or real.

Next we consider collision chains with fixed energy $E$ and angular
momentum $G$. Again we say that a periodic modulo rotation collision
chain $\gamma$ is nondegenerate if the corresponding
$(\tbf,T,\gamma,\Phi)\in\hat\Omega_n$ is a nondegenerate modulo  symmetry
critical point of the functional $\calA^{EG}$ on $\hat\Omega_n$. Thus it has only
degeneracy coming from rotation and time translation invariance.

\begin{thm}\label{thm:EG1}
Let $\gamma$ be a $T$-periodic modulo rotation nondegenerate collision
chain with energy $E$ and angular momentum $G$.  Then  for small
$\mu>0$ there exists a periodic modulo rotation orbit $\gamma_\mu$
of system $(H_\mu)$ with energy $E$ and angular momentum $G$ which
$O(\mu)$-shadows $\gamma$:
$$
\gamma_\mu(t+T_\mu)=e^{i\Phi_\mu}\gamma_\mu(t),\quad \gamma_\mu(t)=\gamma(t)+O(\mu),\qquad t\in[0,T],
$$
where
$$
T_\mu=T+O(\mu),\quad \Phi_\mu=\Phi+O(\mu).
$$
\end{thm}

An estimate (\ref{eq:O(mu)}) holds also here.
Even if $\gamma$ is periodic ($\Phi\in2\pi\Q$), in general the
shadowing orbit $\gamma_\mu$ will be periodic only modulo rotation,
and thus quasiperiodic in a fixed coordinate frame.

To use Theorems \ref{thm:E1}--\ref{thm:EG1}, we need to find nondegenerate modulo symmetry  collision chains.
 In general this is not easy.
A simple application of Theorem \ref{thm:EG1},
based on a perturbative approach, is given in
section \ref{sec:appl}.
More complex applications will be given in a future publication.

In section \ref{sec:action} a description of nondegenerate collision
orbits is given. Using this description, in section \ref{sec:discr}
we  reduce the action functionals  to their discrete versions. Then
in section \ref{sec:proof} we formulate  a local connection result
-- Theorem \ref{thm:connect} -- and use it to prove Theorem
\ref{thm:E1}. The proof of Theorem \ref{thm:EG1} is similar. In
section \ref{sec:Levi} we use Levi-Civita regularization to reduce
Theorem \ref{thm:connect} to Theorem \ref{thm:Shil} which is a
generalization of the Shilnikov Lemma \cite{Turayev} to Hamiltonian
systems with a normally hyperbolic critical manifold.

\begin{rem} In this paper we do not attempt to use global variational methods.
The reason is that although one can use global methods to find
critical points of the action functionals, in general it is hard to
check that the critical points satisfy the changing direction
condition.
\end{rem}

\section{Restricted elliptic limit}

\label{sec:appl}

Suppose that one of the small masses $m_1,m_2$ is much smaller than
the other: $\alpha_1\ll\alpha_2$. In the formal limit $\alpha_1\to 0$ we
obtain the restricted elliptic 3 body problem for which many second species
periodic solutions were obtained in \cite{Bol:ELL}. These results do
not immediately extend to the case of small $\alpha_1>0$. However,
we will show that they can be used to obtain many second species
periodic solutions for the nonrestricted 3 body problem.

Let us fix energy $E$ and angular momentum $G$.
For $\alpha_1\ll \alpha_2$, the Maupertuis--Routh action functional
\begin{equation}\label{eq:pert}
\calA^{EG}(\tbf,T,\gamma,\Phi)=\alpha_1 B(T,\gamma_1)+
\alpha_2 B(T,\gamma_2)+ET-G\Phi,\quad \gamma=(\gamma_1,\gamma_2),
\end{equation}
on $\hat\Omega_n$ is a small perturbation of the Maupertuis--Routh
action functional for the Kepler problem. Indeed,
the functional $\calA_0^{EG}= \calA^{EG}|_{\alpha_1=0}$ does not
depend on $\gamma_1$:
\begin{equation}
\label{eq:B2}
\calA_0^{EG}(\tbf,T,\gamma,\Phi)=\calB^{EG}(T,\gamma_2,\Phi)=B(T,\gamma_2)+TE-G\Phi.
\end{equation}
The condition $\gamma(t_j)\in\Delta$ imposes no
restrictions on $\gamma_2$. Thus the functional $\calB^{EG}$ is
defined on the set $\Pi$ of $(T,\sigma,\Phi)$, where $\sigma:\R\to
U$ is a $W_\loc^{1,2}$ curve such that
$\sigma(t+T)=e^{i\Phi}\sigma(t)$. We have a submersion
$\pi:\hat\Omega_n\to \Pi$, $(\tbf,T,\gamma,\Phi)\to
(T,\gamma_2,\Phi)$, and $\calA_0^{EG}=\calB^{EG}\circ\pi$.

The functional $\calB^{EG}$ is very degenerate, because all orbits
of the Kepler problem with energy $E<0$ are periodic with the same
period $\tau=2\pi (-2E)^{-3/2}$. Suppose $E, G$ are such that there
exists an elliptic  orbit $\Gamma$ of Kepler's problem with energy
$E$ and angular momentum $G$. For definiteness let $G>0$. Then
$0<(-2E)G<1$. The major semiaxis and eccentricity  of $\Gamma$ are
$$
a=(-2E)^{-1},\quad e=\sqrt{1+2EG}.
$$
The Maupertuis action is
$$
J_E(\Gamma)=\int_\Gamma y\cdot dx= 2\pi(-2E)^{-1/2}.
$$
The counterclockwise elliptic orbit $\Gamma:\R\to U$ is defined uniquely
modulo rotation and time translation.

\begin{prop} Let $E<0$, $G>0$ and $(-2E)G<1$.
Then all critical points $\omega=(T,\sigma,\Phi)$ of the functional $\calB^{EG}$ on $\Pi$ belong
to one of the nondegenerate critical manifolds $M_m\subset\Pi$,
$m\in\N$, obtained from $(m\tau,\Gamma,0)$ by rotation and time
translation of $\Gamma$. We have
$$
\calB^{EG}|_{M_m}=B(m\tau,\Gamma)+m\tau E=mJ_E(\Gamma)=2\pi m(-2E)^{-1/2}.
$$
\end{prop}

\proof Let $(T,\sigma,\Phi)\in\Pi$ be a critical point of $\calB^{EG}$.
Then $\sigma$ is a solution of the Kepler problem with energy $E$
and angular momentum $G$ and hence $\sigma$ is a time translation and rotation of
$\Gamma$. Since $\Gamma$ is a non-circular orbit, quasiperiodicity
condition $\sigma(t+T)=e^{i\Phi }\sigma(t)$ implies that
$\Phi=0\bmod 2\pi\Z$ and $T=m\tau$ for some $m\in\N$.

Next we need to check that $M_m$ is a nondegenerate critical
manifold of $\calB^{EG}$. Essentially  this is the same statement,
but now we need to consider the linearized Kepler problem.

The second variation $d^2\calB^{EG}(\omega)$ at $\omega=(m\tau,\Gamma,0)$
is a bilinear form on the tangent space $T_{\omega}\Pi$ which is the set
of $\eta=(\theta,\xi,\phi)$, where $\theta,\phi\in\R$ and
$\xi:\R\to\R^2$ is a vector field such that
\begin{equation}
\label{eq:xi}
\xi(t+m\tau)=\xi(t)+\dot\Gamma(t)\theta+  i\Gamma(t)\phi.
\end{equation}
The standard calculus of variations implies that if $\eta\in T_\omega\Pi$
belongs to the kernel of
$d^2\calB^{EG}(\omega)$, then $\xi$ is a solution of the
variational equation for $\Gamma$ which lie on the zero levels of
the linear first integrals corresponding to the integrals of angular
momentum and energy.
The linear approximations at $\Gamma$ to the integrals of energy and angular momentum are
$$
\dot\Gamma(t)\cdot\dot\xi(t)-\ddot\Gamma(t)\cdot \xi(t)\equiv 0,\quad
i\Gamma(t)\cdot\dot\xi(t)-i\dot\Gamma(t)\cdot \xi(t)\equiv 0.
$$
Condition (\ref{eq:xi}) gives
$$
i\Gamma(t)\cdot \ddot\Gamma(t)\theta\equiv 0,\quad  i\Gamma(t)\cdot \ddot\Gamma(t)\phi\equiv 0.
$$
Since $\Gamma$ is noncircular, $\theta=\phi=0$ and so
$\xi(t+m\tau)=\xi(t)$. It follows that $\eta=(0,\xi,0)$ is
tangent to $M_m$, i.e.\ the variation $\xi(t)$
is obtained by time translation and rotation of
$\Gamma(t)$. \qed

\medskip

The critical manifold $N_m=\pi^{-1}(M_m)\subset\hat\Omega_n$ of
$\calA_0^{EG}$ corresponding to $M_m$ is (up to time translation
and rotation)
$$
N_m=\{(\tbf,m\tau,\sigma,\Gamma,0)\in\hat\Omega_n:
\sigma(t+m\tau)=\sigma(t),\; \sigma(t_j)=\Gamma(t_j)\}.
$$
This is an infinite dimensional nondegenerate critical manifold of
$\calA_0^{EG}$. For nonzero $\alpha_1$, by (\ref{eq:pert}),
\begin{eqnarray*}
\calA^{EG}|_{N_m}=\alpha_1 (B(T,\sigma)+Em\tau - 2\pi
m(-2E)^{-1/2})+ 2\pi m(-2E)^{-1/2}.
\end{eqnarray*}

 By a standard property of nondegenerate critical manifolds
\cite{Palais}, any nondegenerate modulo symmetry critical point
$\omega\in N_m$ of $\calA^{EG}|_{N_m}$ for small $\alpha_1>0$ gives
a nondegenerate modulo symmetry critical point
of $\calA^{EG}$, and hence a nondegenerate modulo
symmetry collision chain with energy $E$ and angular momentum $G$.

Up to an additive constant  and a constant multiple,
$\calA^{EG}|_{N_m}$  is Hamilton's action
$B(m\tau,\sigma)$ for the Kepler
problem. It is defined on the set $\Pi_{\Gamma,m}$
of $(\tbf,\sigma)$, where $\sigma:\R\to
U$ is an $m\tau$-periodic curve such that $\sigma(t_j)=\Gamma(t_j)$.
Thus $B(m\tau,\sigma)=\calB_{\Gamma,m}(\tbf,\sigma)$ is precisely the action
functional whose critical points are collision chains of the
elliptic restricted 3 body problem. This functional was studied in
\cite{Bol:ELL}, and many its nondegenerate critical points were
found for small eccentricity (almost circular  $\Gamma$), i.e.\
$(-2E)G$ close to 1. Also the changing direction and no early
collisions condition  was verified in \cite{Bol:ELL}, and this
carries out for small $\alpha_1>0$. We obtain

\begin{thm} Let
$0<(-2E)|G|<1$ be close to 1. Then for sufficiently small
$\alpha_1>0$ there exist many collision chains $\gamma$ such that
for sufficiently small $\mu>0$, $\gamma$ is $O(\mu)$-shadowed by a
second species periodic modulo rotation solution $\gamma_\mu$ of the
nonrestricted 3 body problem with given $E,G$.
\end{thm}

This result can be improved by using  a more quantitative statement from
\cite{Bol:ELL}. The obtained second species
solutions are periodic in a rotating coordinate frame and
quasiperiodic in a fixed coordinate frame. Proper periodic orbits
will be obtained in a future publication; for them reduction to the
restricted elliptic problem is impossible.

\section{Collision action function}

\label{sec:action}

Collision chains can be represented  as critical points of a
function of a finite number of variables -- discrete action
functional. This is needed for the proof of Theorems \ref{thm:E1}--\ref{thm:EG1}
and in subsequent publications.
Since collision chains are concatenations of  collision
orbits, we need to describe  collision orbits first.

A collision orbit $\gamma=(\gamma_1,\gamma_2)$ of system $(H_0)$ is
a pair of Kepler orbits joining the points $x_-,x_+\in U$. Thus
description of collision orbits is reduced to the classical
Lambert's problem \cite{Whitt}  of joining the points $x_-,x_+$ by a
Kepler orbit.

First  we join the points $x_-,x_+$ by a Kepler orbit
$\Gamma:[0,\tau]\to U$ with fixed energy $E<0$, or, equivalently,
fixed major semi axis $a=(-2E)^{-1}$. Due to scaling invariance of
Kepler's problem without loss of generality set $a=1$. Then a
Kepler ellipse passing through $x_-,x_+$  is determined by the
second focus $F$ such that
\begin{equation}
\label{eq:F} |x_-|+|x_--F|=2,\quad |x_+|+|x_+-F|=2.
\end{equation}
The solution $F=F(x_-,x_+)$ of these equations exists and smoothly
depends on $x_\pm$ if the corresponding circles intersect
transversely, i.e.\ $(x_-,x_+)$ lie in the set
$$
X=\{(x_-,x_+)\in U^2:||x_+|-|x_-||<|x_+-x_-|<4-|x_-|-|x_+|\}.
$$
For $(x_-,x_+)\in X$ there exist two solutions $F$ of equations
(\ref{eq:F}), and we take one of them, for definiteness the one on
the left side of the segment $x_-x_+$.

Let $\Gamma(x_-,x_+)$ be the counter clock wise simple arc of the
constructed Kepler ellipse   joining the points $x_-$ and $x_+$.
Let
$$
f(x_-,x_+)=\int_\Gamma y\cdot dx=\int_\Gamma(2|x|^{-1}-1)^{1/2}|dx|
$$

be the Maupertuis action of $\Gamma$. This is a smooth rotation
invariant function on $X$:
$$
f(e^{i\theta}x_-,e^{i\theta}x_+)=f(x_-,x_+).
$$

\begin{rem} By Lambert's Theorem \cite{Whitt}, $f$ is a function of
$s_\pm=|x_-|+|x_+|\pm|x_--x_+|$ only. An explicit formula  is
\begin{equation}
\label{eq:Lambert}
f(x_-,x_+)=W(s_+)\pm W(s_-),
\end{equation}
where
$$
W(s)=\frac12\sqrt{(4-s)s}+2\arctan\sqrt{\frac{s}{4-s}}.
$$
Plus is taken if $x_+=e^{i\theta}x_-$ with $\theta\in
[\pi,2\pi)$ and minus  if $\theta\in (0,\pi]$. One can check that
$f$ is smooth at $\theta=\pi$.
\end{rem}

Due to scaling invariance of the Kepler problem, for arbitrary negative
energy $E<0$, the Maupertuis action of a simple counter clock wise  arc
$\Gamma=\Gamma(E,x_-,x_+)$ connecting the points $(x_-,x_+)\in X_E=
(-2E)^{-1}X$ is
\begin{eqnarray*}
f(E,x_-,x_+)=\int_\Gamma y\cdot dx =\int_\Gamma
(2(|x|^{-1}+E))^{1/2}|dx|\\= (-2E)^{-1/2}f((-2E)x_-,(-2E)x_+).
\end{eqnarray*}

Next we describe  Kepler orbits $\Gamma:[0,\tau]\to U$
connecting $x_-,x_+$ while making $n=[\Gamma]$  full revolutions around
$0$. To define the number $n\in\Z$ of revolutions, set
$$
n=[(\theta(\tau)-\theta(0))/2\pi],\qquad \Gamma(t)=r(t)e^{i\theta(t)}.
$$

\begin{prop}\label{prop:Lambert_E}
For any $(x_-,x_+)\in X_E$ and any $n\in\Z$:
\begin{itemize}
\item
There exists a Kepler orbit $\Gamma=\Gamma_n(E,x_-,x_+):[0,\tau]\to
U$  of energy $E$ joining the points $x_-,x_+$ and making
$[\Gamma]=n$ revolutions.
\item $\Gamma$ smoothly depends on
$E,x_-,x_+$.
\item
The Maupertuis action of $\Gamma$ is
\begin{eqnarray}
J_n(E,x_-,x_+)=\int_\Gamma y\cdot dx= \nonumber\\
 \label{eq:J_n}
 =(-2E)^{-1/2}(2\pi|n|+(\sgn n)f((-2E)x_-,(-2E)x_+)).
\end{eqnarray}

\item
\begin{equation}
\label{eq:tau1} \tau=\tau_n(E,x_-,x_+)=\frac{\partial}{\partial
E}J_n(E,x_-,x_+).
\end{equation}
\item
\begin{equation}
\label{eq:d2J}
\frac{\partial^2}{\partial
E^2}J_n(E,x_-,x_+)>0,\qquad (x_-,x_+)\in X_E.
\end{equation}
\end{itemize}
\end{prop}

The orbits $\Gamma_n(E,x_-,x_+)$ are nondegenerate (have
non-conjugate end points) and any nondegenerate connecting orbit
with $E<0$ is obtained in this way.

For the classical Lambert's problem \cite{Whitt}, when $\Gamma$ is a simple
elliptic arc, $n=0$ or $n=-1$ depending on if $\Gamma$ is a
counterclockwise or clockwise. We set $\sgn 0=1$.

The first term in (\ref{eq:J_n}) is the Maupertuis action for $n$
complete revolutions around the Kepler ellipse, and the second is
the action of a simple elliptic arc. Equation (\ref{eq:tau1})
follows from the first variation formula; it is essentially Kepler's
time equation. So only inequality (\ref{eq:d2J}) is non-evident. It
is enough to check it for  the classical Lambert's problem with
$n=0,1$. Then (\ref{eq:d2J}) can be deduced from the explicit
formula (\ref{eq:d2F}), although the computation is not trivial. An
equivalent statement was proved in \cite{Simo}. \qed

\medskip

Next we  consider Lambert's problem for fixed time $\tau>0$. This
problem involves solving the transcendental Kepler's equation
so there is no explicit formula for the solution.
Let $D_n\subset \R_+\times U^2$ be the open set which is the image of the diffeomorphism
$$
(E,x_-,x_+)\to (\tau_n(E,x_-,x_+),x_-,x_+),\qquad E<0,\;
(x_-,x_+)\in X_E.
$$

\begin{prop}\label{prop:Lambert1}
For any $n\in\Z$ and any $(\tau,x_-,x_+)\in D_n$:
\begin{itemize}
\item There exists a nondegenerate
Kepler orbit $\Gamma=\Gamma_n(\tau,x_-,x_+):[0,\tau]\to U$ with
$[\Gamma]=n$ full revolutions joining the points $x_-$ and $x_+$.
\item
$\Gamma$ smoothly depends on $\tau,x_-,x_+$.
\item
Hamilton's action
$$
B(\tau,\Gamma)=F_n(\tau,x_-,x_+)
$$
is a smooth function on $D_n$ and
\begin{equation}
\label{eq:d2F}
\frac{\partial^2}{\partial\tau^2}F_n(\tau,x_-,x_+)<0,\qquad
(\tau,x_-,x_+)\in D_n.
\end{equation}
\item All nondegenerate connecting orbits with $E<0$ are $\Gamma_n(\tau,x_-,x_+)$ for some $n\in\Z$
and $(\tau,x_-,x_+)\in D_n$.
\end{itemize}
\end{prop}

\proof We need to find the energy $E<0$ such that the connecting orbit
$\Gamma=\Gamma_n(E,x_-,x_+)$ in Proposition \ref{prop:Lambert_E} has
given time $\tau=\tau_n(E,x_-,x_+)$.  Then
$$
(F_n(\tau,x_-,x_+)+\tau E)\Big|_{\tau=
\tau_n(E,x_-,x_+)}=J_n(E,x_-,x_+).
$$
Hence $J_n$ and $-F_n$ are Legendre transforms of each
other:
\begin{eqnarray*}
J_n(E,x_-,x_+)=\max_{\tau}(F_n(\tau,x_-,x_+)+\tau E),\\
F_n(\tau,x_-,x_+) =\min_E(J_n(E,x_-,x_+)-\tau E).
\end{eqnarray*}
Since $J_n(E,x_-,x_+)$ is convex with respect to $E$,
its Legendre transform $-F_n(\tau,x_-,x_+)$ is convex in $\tau$ and smooth.
\qed

\medskip

The initial and final total momenta of $\gamma$ are given by the
first variation formula
\begin{equation}
\label{eq:DF} y_+=\frac{\partial}{\partial
x_+}F_n(\tau,x_-,x_+),\quad y_-=-\frac{\partial}{\partial
x_-}F_n(\tau,x_-,x_+).
\end{equation}

Now it is easy to describe nondegenerate  collision orbits $\gamma=(\gamma_1,\gamma_2)$ of system $(H_0)$.
Denote by $k=[\gamma]=(k_1,k_2)\in \Z^2$, $k_j=[\gamma_j]$, the rotation vector of $\gamma$.
We obtain

\begin{prop} \label{prop:action1}
For any $k=(k_1,k_2)\in\Z^2$ and any
$(\tau,x_-,x_+)\in V_k=D_{k_1}\cap D_{k_2}$:
\begin{itemize}
\item
There exists a nondegenerate
collision orbit $\gamma:[0,\tau]\to U^2$, $\gamma=\gamma(k,\tau,x_-,x_+)$, with collision points
$\gamma(0)=(x_-,x_-)$, $\gamma(\tau)=(x_+,x_+)$.
\item
$\gamma$ smoothly depends on $(\tau,x_-,x_+)\in V_k$.
\item Hamilton's action of $\gamma$ is
\begin{eqnarray}\label{eq:Sk}
S_k(\tau,x_-,x_+)=A(\tau,\gamma)
=\alpha_1F_{k_1}(\tau,x_-,x_+)+\alpha_2F_{k_2}(\tau,x_-,x_+).
\end{eqnarray}
\item
\begin{equation}
\label{eq:d2S}
\frac{\partial^2}{\partial\tau^2}S_k(\tau,x_-,x_+)<0,\qquad
(\tau,x_-,x_+)\in V_k.
\end{equation}
\end{itemize}
\end{prop}

By the first variation formula \cite{AKN},
$$
dS_k(\tau,x_-,x_+)= y_+\cdot dx_+-y_-\cdot dx_- -
E\,d\tau,
$$
where $E$ is the total energy of the collision orbit $\gamma$, and
$y_\pm=y(t_\pm)$ are total momenta at collisions. Thus
\begin{equation}
\label{eq:gen}
y_+=\frac{\partial S_k}{\partial x_+}(\tau,x_-,x_+),\quad
y_-=-\frac{\partial S_k}{\partial x_-}(\tau,x_-,x_+),\quad E=
-\frac{\partial S_k}{\partial \tau}(\tau,x_-,x_+).
\end{equation}

\begin{rem} We will not need this in the present paper, but for almost all
$(\tau,x_-,x_+)\in V_k$ the collision action
satisfies the twist condition
\begin{equation}\label{eq:twist}
\det\frac{\partial^2 S_k}{\partial x_-\partial x_+}(\tau,x_-,x_+)\ne
0.
\end{equation}
Thus   $S_k$ is the
generating function of a symplectic collision map
$(\tau,x_-,y_-)\to (\tau,x_+,y_+)$.
This will be important for the study of chaotic collision chains.
\end{rem}

\begin{rem}
It can happen that the collision orbit $\gamma$ has early collisions: $\gamma(t)\in\Delta$
for some $t\in(0,\tau)$. To avoid this, we may need to delete from $V_k$ a zero measure set,
see \cite{Bol:ELL}.
\end{rem}

Let us now fix energy $E<0$ and look for collision orbits of system $(H_0)$ with energy
$E$. The map
$$
(\tau,x_-,x_+)\to (-\frac{\partial S_k}{\partial \tau}(\tau,x_-,x_+),x_-,x_+)
$$
is a diffeomorphism of $V_k$ onto an open set $W_k\subset \R\times U^2$.
Let $L_k^E$ be the Legendre transform of $-S_k$ with respect to
$\tau$:
$$
L_k^E(x_-,x_+)=\max_\tau(S_k(\tau,x_-,x_+)+\tau E)=
(S_k(\tau,x_-,x_+)+\tau E)|_{\tau=\tau_k(E,x_-,x_+)},
$$
where $\tau$ is obtained by  solving the last equation (\ref{eq:gen}).
Then $L_k^E$ is a smooth function on
$$
W_k^E=\{(x_-,x_+):(E,x_-,x_+)\in W_k\}.
$$
We obtain

\begin{prop}\label{prop:LE}
Let $E<0$. For any $(x_-,x_+)\in W_k^E$ there exists a unique collision orbit
$\gamma:[0,\tau]\to U^2$
of energy $E$ such that $\gamma(0)=(x_-,x_-)$, $\gamma(\tau)=(x_+,x_+)$.
Its Maupertuis action
$$
L_k^E(x_-,x_+)=\int_\gamma p\cdot dq=\int_\gamma ds_E
$$
is  a smooth function on $W_k^E$. The total momenta at collision are
\begin{equation}
\label{eq:gen_E} y_+=\frac{\partial L_k^E}{\partial
x_+}(x_-,x_+),\quad y_-=-\frac{\partial L_k^E}{\partial
x_-}(x_-,x_+).
\end{equation}
\end{prop}

In terms of actions functions (\ref{eq:J_n}) for the Kepler problem,
$$
L_k^E(x_-,x_+)=\min_{E_1+E_2=E}(\alpha_1
J_{k_1}(E_1,x_-,x_+)+\alpha_2 J_{k_2}(E_2,x_-,x_+)).
$$

\begin{rem} Due to homogeneity of the Kepler problem,
$$
L_k^E(x_-,x_+) =(-2E)^{-1/2}L_k((-2E)x_-,(-2E)x_+),
$$
where $L_k$  corresponds to energy $E=-1/2$.
\end{rem}

\begin{rem} The action functions $S_k$ and $L_k^E$ can not be expressed in elementary functions.
However they admit simple asymptotic representation for large $k$.
This will be done in a subsequent publication.
\end{rem}

In the next section we use the action functions $S_k$ and $L_k^E$ to represent collision
chains as critical points of discrete action functionals.

\section{Discrete variational principles}

\label{sec:discr}

For a given sequence $\kbf=(k^1,\dots,k^n)\in\Z^{2n}$, $k^j\in \Z^2$, define a
{\it discrete Hamilton's action} by
$$
A_\kbf(\sbf,\xbf)=\sum_{j=1}^n S_{k^j}( s_j,x_j,x_{j+1}),
$$
where $\sbf=( s_1,\dots, s_n)\in\R_+^n$, $\xbf=(x_1,\dots ,x_n)\in
U^n$ and $S_{k^j}$ is the action function on $V_{k^j}$ defined by (\ref{eq:Sk}).
The domain of $A_\kbf$ is
$$
V_\kbf=\{(\sbf,\xbf)\in \R_+^n\times U^n: ( s_j,x_j,x_{j+1})\in
V_{k^j},\; x_{n+1}=x_1\}.
$$

Any $(\sbf,\xbf)\in V_\kbf$ defines
$(\tbf,T,\gamma)\in\Omega_n$ as follows. Take $\tbf=(t_1,\dots,t_n)$
so that $s_j=t_{j+1}-t_j$ and set
$$
\gamma(t)=\gamma(k^j,
s_j,x_j,x_{j+1})(t-t_j),\qquad t_j\le t\le t_{j+1},
$$
where $\gamma(k^j,s_j,x_j,x_{j+1}):[0,s_j]\to U^2$ is the
collision orbit in Proposition \ref{prop:action1}. Then
$\gamma=\gamma_\kbf(\sbf,\xbf)$ is a broken trajectory of system
$(H_0)$ with period $T=\sum_{j=1}^n s_j$  and
 Hamilton's action $A_\kbf(\sbf,\xbf)=A(T,\gamma)$. Of
course $(\tbf,\gamma)$ is defined modulo time translation, so we
identify curves which differ by time translation. Thus we defined an
embedding $\iota:V_\kbf\to\widetilde\Omega_n=\Omega_n/\R$  and
$A_\kbf=\calA\circ\iota$.

For collision chains with fixed period $T$, we restrict $A_\kbf$ to
$$
V_\kbf^T=\{(\sbf,\xbf)\in V_\kbf:  \sum_{j=1}^n s_j=T\}.
$$

\begin{prop}\label{prop:T}
Any  critical point  $(\sbf,\xbf)$ of
$A_\kbf$ on $V_\kbf^T$ defines a $T$-periodic collision chain
$\gamma=\gamma_\kbf(\sbf,\xbf)$.
\end{prop}

Indeed, critical points  of $A_\kbf$ satisfy
\begin{eqnarray}
\label{eq:Dx}
\frac{\partial S_{k^j}}{\partial x_j}( s_j,x_j,x_{j+1})+ \frac{\partial S_{k^{j-1}}}{\partial x_j}( s_{j-1},x_{j-1},x_{i})= 0,\\
\label{eq:-E} \frac{\partial S_{k^j}}{\partial
 s_j}( s_j,x_j,x_{j+1})=-E,
\end{eqnarray}
where $-E$ is the Lagrange multiplier. By (\ref{eq:gen}), $E$ is the energy of the
corresponding collision chain $\gamma$. The total momentum at
collision is
$$
y_j=-\frac{\partial S_{k^j}}{\partial x_j}(
s_j,x_j,x_{j+1})=\frac{\partial S_{k^{j-1}}}{\partial x_j}(
s_{j-1},x_{j-1},x_{j}).
$$
Thus $\gamma$ satisfies (\ref{eq:mom})--(\ref{eq:en}).
\qed

\medskip

Proposition \ref{prop:T} follows also from Proposition
\ref{prop:Ham}. Indeed, the functional $A_\kbf$ is the restriction
of Hamilton's action $\calA$ to the set
$\iota(V_\kbf)\subset\widetilde\Omega_n$ of broken extremals. This
set is obtained by equating to zero the differential of $\calA$ for
fixed $\tbf,T,\xbf$.

In Proposition \ref{prop:T} the period $T$ is fixed. For collision
chains with fixed energy $E<0$ we consider a discrete Maupertuis
action functional  on $V_\kbf$:
$$
A_\kbf^E (\sbf,\xbf)= A_\kbf (\sbf,\xbf)+ET,\qquad T=\sum_{j=1}^N
s_j.
$$
Now $T$ is a function on $V_\kbf$, and $A_\kbf^E (\sbf,\xbf)$ is the
Maupertuis action (\ref{eq:AE}) of the broken trajectory
$\gamma=\gamma_\kbf(\sbf,\xbf)$. We obtain

\begin{prop}\label{prop:E}
To any critical point  $(\sbf,\xbf)$ of $A_\kbf^E$ on $V_\kbf$ there
corresponds a periodic collision chain $\gamma$ with  energy $E$.
All nondegenerate  collision chains with energy $E$ are obtained in
this way from nondegenerate modulo rotation critical points of some
$A_\kbf^E$.
\end{prop}

\begin{rem}
Hamilton's action is invariant under rotations:
$A_\kbf(\sbf,e^{i\theta}\xbf)=A_\kbf(\sbf,\xbf)$.
Thus every  critical point of the functional $A_\kbf^E$
 is degenerate.   To obtain nondegenerate critical points we
should consider the quotient functional $\widetilde A_\kbf^E$ on the
quotient space
$$
\widetilde V_\kbf=V_\kbf/\T\subset \R_+^n\times
\widetilde{U^n},\qquad \widetilde{U^n}=U^n/\T\cong \R_+^n\times\T^{n-1}.
$$
\end{rem}

Let us now fix energy $E<0$ and angular momentum $G$ and
consider periodic modulo rotation collision chains $\gamma$ with
given $E,G$. We obtain
the discrete Maupertuis--Routh action functional
$$
A_\kbf^{EG}(\sbf,\xbf,\Phi)=\sum_{j=1}^n S_{k^j}(
s_j,x_j,x_{j+1})+ET-G\Phi,\qquad x_{n+1}=e^{i\Phi}x_1,\quad
T=\sum_{j=1}^N s_j.
$$
The independent variables are $\sbf=( s_1,\dots, s_n)$,
$\xbf=(x_1,\dots,x_n)$ and $\Phi$, so the domain of $A_\kbf^{EG}$ is
$$
\hat V_\kbf=\{(\sbf,\xbf,\Phi)\in\R_+^n\times U^n\times\R: (
s_j,x_j,x_{j+1})\in V_{k^j},\; x_{n+1}=e^{i\Phi}x_1\}.
$$

\begin{prop}\label{prop:Ehat}
To any critical point  $(\sbf,\xbf,\Phi)$ of $A_\kbf^{EG}$ there corresponds
a periodic modulo rotation collision chain $\gamma=\gamma_\kbf(\sbf,\xbf,\Phi)$ with  energy
$E$ and  angular momentum $G$. Any nondegenerate  periodic modulo rotation collision chain
with energy $E$ and angular momentum $G$ is obtained from a nondegenerate modulo rotation
critical point of some $A_\kbf^{EG}$.
\end{prop}

To construct orbits of system $(H_\mu)$ shadowing the collision
chain $\gamma$ corresponding to a critical point $(\sbf,\xbf)$,
we need  to verify the changing direction
condition. For $k=(k_1,k_2)\in \Z^2$ denote
$$
R_k(\tau,x_-,x_+)=F_{k_1}(\tau,x_-,x_+) -F_{k_2}(\tau,x_-,x_+) .
$$
By (\ref{eq:DF}) the relative collision velocities (\ref{eq:v}) of
a collision orbit $\gamma=\gamma(k,\tau,x_-,x_+):[0,\tau]\to U^2$  are
given by
$$
\dot u(0)=-\frac{\partial R_k}{\partial x_-}(\tau,x_-,x_+), \quad
\dot u(\tau)=\frac{\partial R_k}{\partial x_+}(\tau,x_-,x_+).
$$
Thus the changing direction condition for the collision chain
corresponding to $(\sbf,\xbf)$ can be
expressed as follows:
\begin{equation}
\label{eq:change} \frac{\partial R_{k^j}}{\partial x_j}(
s_j,x_j,x_{j+1})+ \frac{\partial R_{k^{j-1}}}{\partial x_j}(
s_{j-1},x_{j-1},x_{j})\ne 0.
\end{equation}

We have
\begin{equation}
\label{eq:AB}
A_\kbf=\alpha_1B_{\kbf_1}+\alpha_2 B_{\kbf_2},
\end{equation}
where $\kbf=(\kbf_1,\kbf_2)$ with
$\kbf_j=(k^1_j,\dots,k^n_j)\in\Z^n$ and
\begin{equation}
\label{eq:Bk} B_{\kbf_j}(\sbf,\xbf)=\sum_{i=1}^n F_{k^i_j}(
s_i,x_i,x_{i+1})
\end{equation}
is the discrete action functional for the Kepler problem.  If
$(\sbf,\xbf)$ is a critical point of $A_\kbf$ with respect to $\xbf$, then by (\ref{eq:Dx}), the
changing direction condition (\ref{eq:change}) is equivalent to
\begin{equation}
\label{eq:change2} \frac{\partial}{\partial
x_j}B_{\kbf_1}(\sbf,\xbf)\ne 0,\qquad
j=1,\dots, n.
\end{equation}

Next we  reformulate the shadowing Theorems \ref{thm:E1}--\ref{thm:EG1}.

\begin{thm}\label{thm:E}
Let $(\sbf,\xbf)\in V_\kbf$ be a  nondegenerate modulo rotation critical point
of $A_\kbf^E$ satisfying the changing direction condition
(\ref{eq:change2}). Then for sufficiently small $\mu>0$ the
corresponding $T$-periodic collision chain $\gamma$ is $O(\mu)$-shadowed modulo
time translation by an almost collision $T_\mu$-periodic orbit
$\gamma_\mu$ of the 3 body problem with period $T_\mu=T+O(\mu)$.
\end{thm}

\begin{thm}\label{thm:EG}
Let $(\sbf,\xbf,\Phi)\in \hat V_\kbf$ be a  nondegenerate modulo rotation critical
point of $A_\kbf^{EG}$ satisfying the changing direction
condition (\ref{eq:change2}). Then for sufficiently small $\mu>0$,
the corresponding collision chain $\gamma$ is $O(\mu)$-shadowed
modulo rotation and time translation by an almost collision periodic modulo rotation orbit
$\gamma_\mu$ of the 3 body problem with energy $E$ and
angular momentum $G$.
\end{thm}

These discrete versions of Theorems \ref{thm:E1}--\ref{thm:EG1} are
most suitable for applications. In  a future publication we will use them in \cite{future} to
find many nontrivial second species solutions.

For a dynamical systems reformulation, it  is convenient to introduce Jacobi's discrete action
functional
$$
J_\kbf^E(\xbf)=\sum_{j=1}^n L_{k_j}^E(x_j,x_{j+1}).
$$
It is defined on
$$
W_\kbf^E=\{\xbf=(x_1,\dots,x_n): (x_j,x_{j+1})\in W_{k^j}^E,\;
x_{n+1}=x_1\}.
$$
Equating to 0 the derivatives of $A_\kbf^{EG}(\sbf,\xbf)$ with respect to $\sbf$,
we obtain:

\begin{prop}\label{prop:Jacobi2}
Any nondegenerate modulo symmetry periodic collision chain with
energy $E$ corresponds to a nondegenerate  modulo rotation critical
point $\xbf$ of some $J_\kbf^E$.
\end{prop}

A critical point $\xbf=(x_1,\dots,x_n)$ of $J_\kbf^E$ is a $n$-periodic trajectory of a
discrete Lagrangian system $(\calL^E)$ with a multivalued discrete Lagrangian
$\calL^E=\{L_k^E\}_{k\in\Z^2}$:
$$
 \frac{\partial}{\partial x_j}(L_{k^{j-1}}^E(x_{j-1},x_{j})+L_{k^j}^E(x_j,x_{j+1}))=0,\qquad j=1,\dots,n.
 $$
Thus  description of second species solutions is reduced to the
dynamics of a discrete Lagrangian system $(\calL^E)$.  Under a twist
condition, a periodic trajectory of system $(\calL^E)$ corresponds
to a periodic trajectory $(x_j,y_j)$, $y_j=
-\frac{\partial}{\partial x_j}L_{k^j}^E(x_j,x_{j+1})$, of a sequence
of symplectic twist maps $(x_j,y_j)\to (x_{j+1},y_{j+1})$ with
generating functions $L_{k^j}^E$. We postpone this reformulation to
a future paper, where we deal with chaotic almost collision orbits.

The minimal degeneracy of a critical point of $J_\kbf^E$ is at least 1 due to
rotational symmetry $J_\kbf^E(e^{i\theta}\xbf)= J_\kbf^E(\xbf)$.
This implies that the discrete Lagrangian system (or the corresponding symplectic map)
has an integral of angular momentum
$G=ix_j\cdot y_j$.
One can perform Routh's reduction in this discrete Lagrangian system
reducing it to one degree of freedom \cite{Hill}, but this complicates
the discrete Lagrangian.

For periodic modulo rotation collision chains with fixed $E,G$ we have:

\begin{prop}\label{prop:Jacobi-Routh}
Any nondegenerate periodic modulo rotation collision chain with
energy $E$ corresponds to a nondegenerate  modulo rotation critical
point $(\xbf,\Phi)$ of the discrete Jacobi--Routh action functional
$$
J_\kbf^{EG}(\xbf,\Phi)=J_\kbf^E(\xbf)-G\Phi,\qquad x_{n+1}=e^{i\Phi}x_1.
$$
\end{prop}

The proofs of Theorems \ref{thm:E} and \ref{thm:EG} are
modifications of the proof of Theorem 2.1 in \cite{Bol:shadow}.
They are based on the Levi-Civita regularization and shadowing. The
proof of Theorem \ref{thm:E} will be given in the next section. The
proof of  Theorem \ref{thm:EG} is similar and will be omitted.

\section{Proof of Theorem \ref{thm:E}}

\label{sec:proof}

For $\mu\ne 0$, the action functional (\ref{eq:Amu}) of system
$(H_\mu)$ is singular when $\gamma$ approaches $\Delta$. We will
formulate a variational problem for almost collision orbits of
system $(H_\mu)$ with given energy $E$ which  has
no singularity at $\Delta$.

Let us fix energy $E<0$. Trajectories of system $(H_\mu)$ with
energy $E$ are extremals of the Jacobi action functional
$$
J_\mu^E(\gamma)=\int_\gamma ds_\mu^E,\qquad ds_\mu^E=\max_p\{p\cdot
dq: H_\mu(q,p)=E\}.
$$
Away from $\Delta$, the functional $J_\mu^E$ is a regular perturbation of
the Jacobi functional $J^E$ for system $(H_0)$.
Regularizing $J_\mu^E$ near $\Delta$ requires some preparation.

First we describe local behavior of trajectories of system $(H_0)$
colliding with $\Delta$. We will  use  the variables (this is a
version of Jacobi's variables)
\begin{equation}
\label{eq:xu}
x=\alpha_1q_1+\alpha_2q_2,\quad y=p_1+p_2,\quad u=q_2-q_1,\quad v=\alpha_1p_2-\alpha_2p_1.
\end{equation}
Thus $x$ is the center of mass of  $m_1,m_2$, $y$ is their total
momentum, $u$ is their relative position, and $v$  is the scaled
relative velocity.  The change is symplectic:
$$
p\cdot dq=p_1\cdot dq_1+p_2\cdot dq_2=y\cdot dx+v\cdot du.
$$
The inverse change is
$$
q_1=x- \alpha_2u,\quad q_1=x+ \alpha_1u,\quad p_1= \alpha_1
y-v,\quad p_2= \alpha_2y+v.
$$
For solutions of system $(H_0)$, $y=\dot x$ and $v=\alpha \dot u$,
where $\alpha=\alpha_1\alpha_2$.

Let $\gamma$  be a trajectory of system $(H_0)$ with energy $E$.
We denote by $(x(t),y(t),u(t),v(t))$ its representation
in Jacobi's variables. If $\gamma$  has a collision at $t=0$, i.e.\
$u(0)=0$, $x(0)=x_0$, then
$$
H_0=|y_0|^2/2+|v_0|^2/2\alpha-|x_0|^{-1}=E.
$$
We assume that collisions occurs with nonzero relative speed $v_0\ne 0$.
Then there exists $\delta>0$ such that $(x_0,y_0) $ lies
in a compact set
\begin{equation}
\label{eq:ME} M=M_\delta^E=\{(x_0,y_0): \lambda(x_0,y_0)= E-
|y_0|^2/2+|x_0|^{-1}\ge \delta,\; |x_0|\ge\delta\}.
\end{equation}
We fix $\delta>0$. Eventually it will taken sufficiently small.
Denote $B_\rho=\{u\in\R^2:|u|\le\rho\}$ and $S_\rho=\partial B_\rho$.
We have

\begin{lem}\label{lem:shoot}
Take any $\delta>0$ and let $\rho>0$ be sufficiently small. Then for
any $(x_0,y_0)\in M$ and any $u_+\in S_\rho$ there exists a
trajectory $\gamma_+:[0,\tau_+]\to U^2$  of system $(H_0)$ with
energy $E$ such that:
\begin{itemize}
\item  $u(t)\in B_\rho$ for $0\le t\le \tau_+$
and  $x(0)=x_0$, $y(0)=y_0$, $u(0)=0$,
$u(\tau_+)=u_+$.
\item
 $\gamma_+$  smoothly depends
on $(x_0,y_0,u_+)\in M\times S_\rho$ and
\begin{eqnarray}
\tau_+=\tau_+(x_0,y_0,u_+)=\rho\sqrt{\alpha/2\lambda(x_0,y_0)}+O(\rho^2),\nonumber\\
x(\tau_+)=\xi_+(x_0,y_0,u_+)=x_0+ \rho\sqrt{\alpha/2\lambda(x_0,y_0)}y_0+O(\rho^2)
\label{eq:xi+}.
\end{eqnarray}
\item
The Maupertuis action of $\gamma_+$  has the form
\begin{eqnarray}
\label{eq:a+}
J^E(\gamma_+)=\int_{\gamma_+} p\cdot dq=a_+(x_0,y_0,u_+)\\
=\rho\sqrt{2\alpha/\lambda(x_0,y_0)}(E+|x_0|^{-1})+O(\rho^2).\nonumber
\label{eq:f_E}
\end{eqnarray}
\end{itemize}
\end{lem}

\begin{rem} On $S_\rho$ we use the polar coordinate $\theta$, where $u=\rho e^{i\theta}$.
Thus $O(\rho^2)$ means a function of $x_0,y_0,\theta$ whose  $C^2$
norm  is bounded by $c\rho^2$ with $c$
independent of $\rho$.
\end{rem}

The proof is obtained by a simple shooting argument, because $H_0$ has no singularity at $\Delta$:
$$
x(t)=x_0+ty_0+O(t^2),\quad u(t)=tv_0/\alpha+O(t^2).
$$
It remains to solve the equation $u(\tau_+)=u_+$ for $\tau_+$ and $v_0$,
where $|v_0|=\sqrt{2\alpha\lambda(x_0,y_0)}$. \qed

\medskip

Similarly, we have a trajectory $\gamma_-:[\tau_-,0]\to U^2$ of
system $(H_0)$ with energy $E$ such that $x(0)=x_0$, $y(0)=y_0$,
$u(0)=0$, $u(\tau_-)=u_-$. Then
$$
\tau_-=\tau_-(x_0,y_0,u_-),\quad x(\tau_-)=\xi_-(x_0,y_0,u_-),\quad
J^E(\gamma_-)=a_-(x_0,y_0,u_-).
$$

If $(x_0,y_0)\in M$, then $x_0$ belongs to
$$
D=D_\delta^E=\{x:\delta\le |x|\le (\delta-E)^{-1}\}.
$$
Let  $\Sigma_\rho$ be the boundary of the tubular neighborhood
$N_\rho$ of $D^2\subset\Delta$:
$$
\Sigma_\rho=\{q: x\in D,\; u\in S_\rho\},\quad N_\rho=\{q: x\in D,\; u\in B_\rho\}.
$$
Fix arbitrary large $C>0$ and let $$ K_\rho=\{(q_0,q_+)\in
D^2\times\Sigma_\rho:|q_0-q_+|\le C\rho\}.
$$

\begin{lem}
\label{lem:shoot2} If $\rho>0$ is sufficiently small, then for any
$(q_0,q_+)\in K_\rho$, there exists a trajectory
$\gamma:[0,\tau_+]\to N_\rho$ of system $(H_0)$ with energy $E$
joining $q_0$ with $q_+$. Moreover $\gamma$ smoothly depends on
$(q_0,q_+)$ and its Maupertuis action has the form
\begin{eqnarray}
J^E(\gamma)=d_E(q_0,q_+)\nonumber\\
=\sqrt{2\alpha^{-1}(E+|x_0|^{-1})(|x_+-x_0|^2+\alpha^2\rho^2)}
+O(\rho^2).\label{eq:d_E}
\end{eqnarray}
\end{lem}

Here $d_E(q_0,q_+)$ is the distance in the Jacobi metric $ds_E$.

\medskip

\noindent
{\it Proof.}
The condition $\gamma(\tau_+)=q_+=(x_+,u_+)$ gives
$u(\tau_+)=u_+$, $x(\tau_+)=x_+$. Then (\ref{eq:xi+})
makes it possible to determine
\begin{equation}
\label{eq:eta+}
y_0=\eta_+(q_0,q_+)=\sqrt{\frac{2(E+|x_0|^{-1})}{|x_+-x_0|^2+\alpha^2\rho^2}}
(x_+-x_0)+O(\rho).
\end{equation}
\qed

\medskip

Next we connect points $q_-,q_+\in\Sigma_\rho$ by a reflection trajectory of energy $E$.
Let  $P_\rho=\{(q_-,q_+)\in\Sigma_\rho^2:|q_--q_+|\le C\rho\}$.

\begin{prop}\label{prop:reflect}
Let $\rho>0$ be sufficiently small.
Then for any $(q_-,q_+)\in P_\rho$:
\begin{itemize}
\item
There exist $\tau_-<0<\tau_+$ and a broken trajectory
$\gamma:[\tau_-,\tau_+]\to N_\rho$ with energy $E$ such that
$\gamma(0)=q_0=(x_0,x_0)\in D^2$, $\gamma|_{[0,\tau_+]}$,
$\gamma|_{[\tau_-,0]}$ are trajectories of system $(H_0)$,
$\gamma(\tau_\pm)=q_\pm$ and there is no jump of total momentum at
collision: $y(+0)=y(-0)=y_0$.
\item $\gamma$ smoothly depends on $(q_-,q_+)\in P_\rho$.
\item The Maupertuis action of $\gamma$ has the form
\begin{eqnarray}
J_E(\gamma)=\int_\gamma ds_E=g_E(q_-,q_+)=d_E(q_0,q_-)+d_E(q_0,q_+)\\
=\sqrt{2(E+|x_0|^{-1})(|x_+-x_-|^2+4\alpha^2\rho^2)}
+O(\rho^2).
\end{eqnarray}
\item
\begin{eqnarray}
x_0=\xi(q_-,q_+)=(x_++x_-)/2+O(\rho^2),\label{eq:xi2}\\
y_0=\eta(q_-,q_+)=\sqrt{\frac{2(E+|x_0|^{-1})}{|x_+-x_-|^2+4\alpha^2\rho^2}}
(x_+-x_-)+O(\rho),\nonumber\\
\tau_\pm=\tau_\pm(q_-,q_+)=\pm\sqrt{\frac{|x_+-x_-|^2+4\alpha^2\rho^2}
{8(E+|x_0|^{-1})}}+O(\rho^2).\nonumber
\end{eqnarray}
\end{itemize}
\end{prop}

\proof We find $x_0$ from the equation
$$
\frac{\partial}{\partial x_0}(d_E(q_0,q_-)+d_E(q_0,q_+))=0\;
\Leftrightarrow\;\eta_+(q_0,q_+)=\eta_-(q_0,q_-),
$$
where $\eta_\pm$ is defined in (\ref{eq:eta+}). Differentiating
(\ref{eq:d_E}), we  see that the Hessian matrix
$$
\sqrt{\frac{8(E+|x_0|^{-1})}{\alpha(|x_+-x_-|^2+4\alpha^2\rho^2)}}
\left(I-\frac{(x_+-x_-)\otimes
(x_+-x_-)}{|x_+-x_-|^2+4\alpha^2\rho^2}+O(\rho)\right)
$$
is nondegenerate. By the implicit function theorem, the solution $x_0=\xi(q_-,q_+)$ is smooth.
\qed

\medskip

A  similar result holds for the perturbed system $(H_\mu)$, but it
is no longer easy to prove. Fix  an arbitrary small constant
$\delta>0$.

\begin{lem}\label{lem:connect}
Let $\rho>0$ be sufficiently small. There exists $\mu_0>0$ such that
for all $\mu\in(0,\mu_0]$, any $(x_0,y_0)\in M$ and any $u_\pm\in
S_\rho$ such that  $|u_++u_-|\ge \delta \rho$:
\begin{itemize}
\item
There exist $t_-<0<t_+$ and a trajectory $\gamma:[t_-,t_+]\to
N_\rho$ of system $(H_\mu)$ with energy $E$  such that
$u(t_\pm)=u_\pm$, $x(0)=x_0$, $y(0)=y_0$.
\item
$\gamma$ smoothly depends on $(x_0,y_0,u_-,u_+,\mu)\in M\times
S_\rho^2\times(0,\mu_0]$ and converges to a concatenation of
trajectories $\gamma_\pm$ in Lemma \ref{lem:shoot} as $\mu\to 0$.
\item
The Maupertuis action of $\gamma$ has the form
\begin{eqnarray}
J_\mu^E(\gamma)=\int_\gamma p\cdot dq
=f_\mu^E(x_0,y_0,u_-,u_+)\\
=a_+(x_0,y_0,u_+)+a_-(x_0,y_0,u_-)+\mu \hat
a(x_0,y_0,u_-,u_+,\mu),
\end{eqnarray}
where $\hat a$ is $C^2$ bounded on $M\times
S_\rho^2\times(0,\mu_0]$.
\item
\begin{eqnarray}
t_\pm=\tau_\pm(x_0,y_0,u_\pm)+\mu \hat
\tau_\pm(x_0,y_0,u_+,u_-,\mu),\nonumber\\
x(t_\pm)=x_\mu^\pm(x_0,y_0,u_-,u_+)\nonumber\\
=\xi_\pm(x_0,y_0,u_\pm)+\mu\hat x_\pm(x_0,y_0,u_-,u_+,\mu).
\label{eq:xpm}
\end{eqnarray}
where $\hat \tau_\pm$ and $\hat x_\pm$ are uniformly $C^1$ bounded on
$M\times S_\rho^2\times(0,\mu_0]$.
\item
If $|u_+-u_-|\ge \delta \rho$, then
\begin{equation}
\label{eq:min} \mu a\le \min_{t\in [t_-,t_+]}d(\gamma(t),\Delta)\le
\mu b, \qquad 0<a<b.
\end{equation}
\end{itemize}
\end{lem}

The proof of Lemma \ref{lem:connect} is given in section
\ref{sec:Levi}. It is based on Levi-Civita regularization and a
generalization of Shilnikov's Lemma \cite{Shilnikov}, see also \cite{Turayev}, to
normally hyperbolic critical manifolds of a Hamiltonian system.

Next we deduce a local connection theorem. Fix arbitrary small
$\delta>0$, arbitrary large $C>0$ and let
$$
Q_\rho=\{(q_-,q_+)\in\Sigma_\rho^2:|q_--q_+|\le C\rho,\;
|u_++u_-|\ge \delta \rho\}.
$$

\begin{thm}\label{thm:connect} Let $\rho>0$ be sufficiently small.
There exists $\mu_0>0$ such that for all $(q_-,q_+,\mu)\in Q_\rho\times (0,\mu_0]$:
\begin{itemize}
\item
There exist $t_-<0<t_+$ and a trajectory $\gamma:[t_-,t_+]\to
N_\rho$ of system $(H_\mu)$ with energy $E$ such that
$\gamma(t_\pm)=q_\pm$ and the minimum of $d(\gamma(t),\Delta)$ is
attained at $t=0$.
\item
$\gamma$ smoothly depends on $(q_-,q_+,\mu)\in Q_\rho\times
(0,\mu_0]$ and converges to a reflection trajectory   in Proposition
\ref{prop:reflect} as $\mu\to 0$.
\item
The Maupertuis action of $\gamma$ has the form
$$
J_\mu^E(\gamma)=\int_\gamma p\cdot dq
=g_\mu^E(q_-,q_+)=g_E(q_-,q_+)+\mu \hat g(q_-,q_+,\mu),
$$
where $\hat g$ is $C^2$ bounded on $ Q_\rho\times(0,\mu_0]$.
\item
If $|u_+-u_-|\ge \delta \rho$, then (\ref{eq:min}) holds.
\end{itemize}
\end{thm}

Thus the action $J_\mu^E(\gamma)=g_\mu^E(q_-,q_+)$ has a limit
$g_E(q_-,q_+)$ as $\mu\to 0$ which is the action of the reflection orbit
in Proposition \ref{prop:reflect}.
The condition that the distance to
$\Delta$ is attained at $t=0$ is needed only to exclude time
translations, so that $t_\pm$ are uniquely defined.

\medskip

\noindent {\it Proof.} We need to find $(x_0,y_0)$ such that
$x_\mu^\pm(x_0,y_0,u_-,u_+)=x_\pm$. Since the implicit function theorem worked
in the proof of Proposition \ref{prop:reflect}, by (\ref{eq:xpm}),
for small $\mu>0$ it will work also here.
\qed

\medskip

\noindent {\it Proof of Theorem \ref{thm:E}.} Let $\gamma$ be a
nondegenerate $n$-collision chain with energy $E$. Let
$\tbf=(t_1,\dots,t_n)$ be collision times,
$\xbf=(x_1,\dots,x_n)$, $\gamma(t_j)=(x_j,x_j)$, the corresponding
collision points and $y_j$ the collision total momenta.
Take $\delta>0$ so small that the collision points and the collision speeds
$v_j^\pm=v(t_j\pm 0)$ satisfy
$$
|x_j|\ge\delta,\quad |v_j^+|=|v_j^-|\ge\delta\sqrt{2\alpha}.
$$
Then $(x_j,y_j)\in M$ and $x_j\in D$.

Take small $\rho>0$ and let $t_j^\pm=t_j\pm s_j^\pm$ be the closest to
$t_j$ times when $q_j^\pm=\gamma(t_j^\pm)\in\Sigma_\rho$. Since
$\gamma$ satisfies the changing direction condition,
$(q_j^-,q_j^+)\in Q_\rho$ if $C>0$ is taken sufficiently large and
$\rho>0$ sufficiently small. Moreover  for $\xi_j^\pm\in\Sigma_\rho$
 close to $q_j^\pm$, we have  $(\xi_j^-,\xi_j^+)\in Q_\rho$.
 Thus by Theorem \ref{thm:connect}
 for small $\mu>0$ the points $\xi_j^\pm$ can be joined in $N_\rho$ by a
 trajectory $\gamma_\mu^j$ of system $(H_\mu)$ with energy $E$
 and the Maupertuis action $J_\mu^E(\gamma_\mu^j)=g_\mu^E(\xi_j^-,\xi_j^+)$.

 Since $\gamma|_{[t_j,t_{j+1}]}$ is nondegenerate and, by no early collisions condition,
 does not come near $\Delta$,
 for $\xi_j^+$ close to $q_j^+$ and $\xi_{j+1}^-$ close to  $q_{j+1}^-$ and small $\mu>0$,
 the points  $\xi_j^+$  and $\xi_{j+1}^-$   can be joined by
 a trajectory $\sigma_\mu^j$ of system $(H_\mu)$ with energy $E$ and the Maupertuis action
 $J_\mu^E(\sigma_\mu^E)=h_\mu^E(\xi_j^+,\xi_{j+1}^-)$.
This trajectory smoothly depends on $\mu$ also for $\mu=0$, and
$h_0^E(\xi_j^+,\xi_{j+1}^-)$ is the Maupertuis action of a
connecting trajectory of system $(H_0)$.

 Combine the trajectories $\gamma_\mu^j$, $\sigma_\mu^j$ in
 a broken trajectory $\gamma_\mu$ with energy $E$ and Maupertuis action
 \begin{eqnarray*}
 f_\mu(\xi)=J_\mu^E(\gamma_\mu)=
 \sum_{j=1}^n(J_\mu^E(\gamma_\mu^j)+J_\mu^E(\sigma_\mu^j))\\
 =\sum_{j=1}^n(g_\mu^E(\xi_j^+,\xi_{j+1}^-)+h_\mu^E(\xi_j^-,\xi_j^+)),
 \qquad \xi=(\xi_1^-,\xi_1^+,\dots, \xi_n^-,\xi_n^+)\in \Sigma_\rho^{2n}.
 \end{eqnarray*}

 The function $f_\mu$  has a limit $f_0$ as $\mu\to 0$ and
 $$
 f_\mu(\xi)=f_0(\xi)+\mu\hat f(\xi,\mu),
 $$
where
$$
f_0(\xi)=\sum_{j=1}^n(d(\xi_j,\xi_j^+)+h_0^E(\xi_j^+,\xi_{j+1}^-)+d(\xi_{j+1}^-,\xi_{j+1})),
$$
and $\xi_j=\xi(\xi_j^-, \xi_j^+)$ is defined by (\ref{eq:xi2}). The
remainder $\hat f(\xi,\mu)$ is $C^2$ bounded on $Y\times (0,\mu_0]$,
where the neighborhood $Y\subset \Sigma_\rho^{2n}$  of
$$
 \qbf=(q_1^-,q_1^+,\dots, q_n^-,q_n^+)
 $$
 is independent of $\mu$.

 Looking for critical points with respect to $\xi_j^\pm$ with fixed $\xi_j$ we
obtain $f_0=J_\kbf^E(\xi_1,\dots,\xi_n)$ for some $\kbf\in\Z^{2n}$,
and $J_\kbf^E$ has a nondegenerate modulo rotation critical point
$\xbf$.  Thus $f_0$ has a nondegenerate critical point $\qbf\in
\Sigma_\rho^{2n}$. Then for small $\mu>0$
 the function $f_\mu(\xi)$ has a nondegenerate modulo rotation
 critical point $\xi_\mu$ close to $\qbf$. The corresponding broken trajectory
 $\gamma_\mu$ has no break of velocity at intersection points $\xi_j^\pm$
 with $\Sigma_\rho$ and hence $\gamma_\mu$ is a periodic trajectory of system $(H_\mu)$
 with energy $E$.
 \qed

\section{Levi-Civita regularization}

\label{sec:Levi}

In this section we prove Lemma \ref{lem:connect}.
In the Jacobi variables   (\ref{eq:xu}), the Hamiltonian $H_\mu$ takes the form
$$
H_\mu=\frac{(1+\mu)|y|^2}2+\frac{|v|^2}{2 \alpha}-\frac{ \alpha_1}{|
\alpha_2u-x|} -\frac{ \alpha_2}{| \alpha_1u+x|}-\frac{\mu
\alpha}{|u|}.
$$

Let us perform the Levi-Civita regularization on the fixed energy
level $H_\mu=E$. We identify $u,v\in\R^2=\C$ with
complex numbers and  make a change of
variables
$$
u=\xi^2,\quad v=\eta/2\bar\xi.
$$
Since
$$
v\cdot du=\Real(v\, d\bar u)= \Real(\eta\, d\bar
\xi)=\eta\cdot d\xi,
$$
the change is symplectic:
\begin{equation}
\label{eq:sympl}
 p\cdot dq=y\cdot dx+\eta\cdot d\xi.
\end{equation}

Finally, we obtain a transformation
$$
q_1=x- \alpha_2\xi^2,\quad q_1=x+ \alpha_1\xi^2,\quad p_1= \alpha_1 y-\eta/2\bar\xi,\quad p_2=
\alpha_2y+\eta/2\bar\xi.
$$
The Levi-Civita map
$$
g:\R^2\times\R^2\times U \times\R^2\to
(\R^4\setminus\Delta)\times\R^4,\qquad
g(x,y,\xi,\eta)=(q_1,q_2,p_1,p_2),
$$
is a symplectic double covering
undefined at $\xi=0$ which corresponds to the collision set
$\Delta$.

Let
\begin{eqnarray*}
\calH_\mu^E(x,y,\xi,\eta)=|\xi|^2(H_\mu\circ g-E)+\mu\alpha\\
=\frac{|\eta|^2}{8
\alpha}-|\xi|^2\left( E+\frac{ \alpha_1}{|
\alpha_2\xi^2-x|} +\frac{ \alpha_2}{|
\alpha_1\xi^2+x|}-\frac{(1+\mu)|y|^2}2 \right).
\end{eqnarray*}

Let  $\Sigma_\mu^E=\{H_\mu=E\}$ and
$\Gamma_\mu^E=\{\calH_\mu^E=\mu\alpha\}$. Since
$g(\Gamma_\mu^E)=\Sigma_\mu^E$, the map $g$ takes orbits of system
$(\calH_\mu^E)$  on $\Gamma_\mu^E$ to orbits of system $(H_\mu)$ on
$\Sigma_\mu^E$. The time parametrization is changed: the new time is
given by $d\tau=|\xi|^2dt$. In the following we will continue to
denote the new time by $t$.

The singularity at $\Delta$ disappeared after regularization. The
regularized Hamiltonian $\calH_\mu^E$ is smooth on
$$
\calP=\{(x,y,\xi,\eta)\in U\times\R^6:x\ne \alpha_2\xi^2,\;x\ne-\alpha_1\xi^2\}.
$$
which means excluding collisions of $m_1$ and $m_2$ with $m_3$. The
parameter $\mu\alpha$ may be regarded as new energy. The rotation group
and the integral of angular momentum are now
$$
(x,y,\xi,\eta)\to (e^{i\theta}x,e^{i\theta}
y,e^{i\theta/2}\xi,e^{i\theta/2}\eta),\quad G=ix\cdot y+i\xi\cdot
\eta/2.
$$

The Hamiltonian $\calH_\mu^E$  has a critical manifold $\xi=\eta=0$
which is contained in  the level set   $\Gamma_0^E$ of $\calH_0^E$.
We have
$$
\calH_0^E(x,y,\xi,\eta)=\frac{|\eta|^2}{8
\alpha}-|\xi|^2\lambda(x,y)+O(|\xi|^4).
$$
Collisions of $m_1,m_2$ with nonzero relative velocity correspond to
the solutions asymptotic to
$$
\calM= M\times\{(0,0)\},
$$
where $M$ is as in (\ref{eq:ME}). This is is a compact normally
hyperbolic symplectic critical manifold for $\calH_0^E$. We obtain

\begin{thm}
Collision orbits of system $(H_0)$ with energy $E$ correspond
to orbits of system $(\calH_0^E)$ doubly asymptotic  to $\calM$.
Orbits of system $(H_\mu)$  with energy $E$
passing $O(\mu)$-close to the singular set $\Delta$  correspond to
orbits of system $(\calH_\mu^E)$ on the level $\Gamma_\mu^E$ passing
$O(\sqrt{\mu})$-close to  $\calM$.
\end{thm}

Next we translate Lemma \ref{lem:connect} to the new variables.

Let $r>0$ and let  $\calN_r=M\times
B_r$ be a tubular neighborhood of $\calM$ in $\calP$. By the stable and unstable
manifold theorems for normally hyperbolic invariant manifolds \cite{Fenichel},
if $r>0$ is small enough, for
any $(x_0,y_0)\in M$ and $\xi_-\in S_r$ there exists a unique
solution $\zeta_-:[0,+\infty)\to \calN_r$,
$\zeta_-(t)=(x(t),y(t),\xi(t),\eta(t))$, of system $(\calH_0^E)$ such
that $\xi(0)=\xi_-$ and $\zeta(\infty)=(x_0,y_0,0,0)\in\calM$. We
denote its action by
$$
J(\zeta_-)=\int_{\zeta_-}y\cdot dx+\eta\cdot
d\xi=J_-(x_0,y_0,\xi_-).
$$
Since the stable and unstable manifolds are smooth, $J_-$ is a
smooth function on $M\times S_r$. Similarly we define the  function
$J_+$ on $M\times S_r$ as the action of a solution $\zeta_+$
asymptotic to $\calM$ as $t\to-\infty$.

We have an analog of Shilnikov's Lemma \cite{Shilnikov}. Fix small $\eps>0$ and denote
$$
\calQ_r=\{(x_0,y_0,\xi_-,\xi_+)\in M\times S_r^2: \xi_-\cdot\xi_+\ge \eps^2 r^2\}.
$$

\begin{thm}
\label{thm:Shil} There exists $r>0$ and $\mu_0>0$ such that for
any $(x_0,y_0,\xi_-,\xi_+,\mu)\in \calQ_r\times (0,\mu_0]$:
\begin{itemize}
\item
There exists $T>0$ and a
solution
$$
\zeta(t)=(x(t),y(t),
 \xi(t),\eta(t))\in \calN_r,\qquad t \in [- T, T],
$$
of system $(\calH_\mu^E)$  on $\Gamma_\mu^E$  such that
\begin{equation}
\label{eq:BC0} x(0) = x_0,\quad y(0)=y_0,\quad \xi(- T)=\xi_-,\quad
\xi( T)= \xi_+.
\end{equation}
\item
$\zeta$ smoothly depends on $(x_0,y_0,\xi_-,\xi_+,\mu)\in \calQ_r\times (0,\mu_0]$.
\item
The Maupertuis action is a smooth function on $\calQ_r\times(0,\mu_0]$ and has the form
\begin{equation}
J(\zeta)=\int_\zeta y\cdot dx+\eta\cdot
d\xi=J_-(x_0,y_0,\xi_-)+J_+(x_0,y_0,\xi_+)+\mu\hat
J(x_0,y_0,\xi_-,\xi_+,\mu),
\end{equation}
where $\hat J$ is $C^2$ bounded on $\calQ_r\times (0,\mu_0]$.
\end{itemize}
\end{thm}

A result very similar to Theorem \ref{thm:Shil} was proved in \cite{Bol:shadow}.   A complete proof
of Theorem \ref{thm:Shil} will be published in \cite{future}.

\medskip

\noindent{\it Proof of Lemma \ref{lem:connect}.} We set $\rho=r^2$
and $u=\xi^2$. For given $u_\pm\in S_\rho$ take $\xi_\pm\in S_r$
such that $\xi_+\cdot\xi_-\ge 0$. If $\xi_+\cdot\xi_->0$, then
$u_+\ne -u_-$. Moreover for given $\delta>0$ there exists $\eps>0$
such that $|u_-+u_+|\ge\delta\rho$ implies $\xi_-\cdot\xi_+\ge
\eps^2 r^2$. If $\zeta$ is a trajectory in Theorem \ref{thm:Shil},
then the corresponding trajectory $g(\zeta)$ of system $(H_\mu)$
satisfies the conditions of Lemma \ref{lem:connect}. In particular,
by (\ref{eq:sympl}), $$ J_\pm(x_0,y_0,\xi_\pm)=a_\pm(x_0,y_0,u_\pm).
$$
\qed

\end{document}